

\frenchspacing

\parindent15pt

\abovedisplayskip4pt plus2pt
\belowdisplayskip4pt plus2pt 
\abovedisplayshortskip2pt plus2pt 
\belowdisplayshortskip2pt plus2pt  

\font\twbf=cmbx10 at12pt
 at12pt
 at12pt

\font\sc=cmcsc10

\font\ninerm=cmr9 
\font\nineit=cmti9 
\font\ninesy=cmsy9 
\font\ninei=cmmi9 
\font\ninebf=cmbx9 

\font\sevenrm=cmr7  
 
\font\seveni=cmmi7  
\font\sevensy=cmsy7 

\font\fivenrm=cmr5  
\font\fiveni=cmmi5  
\font\fivensy=cmsy5 

\def\nine{%
\textfont0=\ninerm \scriptfont0=\sevenrm \scriptscriptfont0=\fivenrm
\textfont1=\ninei \scriptfont1=\seveni \scriptscriptfont1=\fiveni
\textfont2=\ninesy \scriptfont2=\sevensy \scriptscriptfont2=\fivensy
\textfont3=\tenex \scriptfont3=\tenex \scriptscriptfont3=\tenex
\def\rm{\fam0\ninerm}%
\textfont\itfam=\nineit    
\def\it{\fam\itfam\nineit}%
\textfont\bffam=\ninebf 
\def\bf{\fam\bffam\ninebf}%
\normalbaselineskip=11pt
\setbox\strutbox=\hbox{\vrule height8pt depth3pt width0pt}%
\normalbaselines\rm}

\hsize30cc
\vsize44cc
\nopagenumbers

\def\luz#1{\luzno#1?}
\def\luzno#1{\ifx#1?\let\next=\relax\yyy
\else \let\next=\luzno#1\xxx\fi\next}
\def\sp#1{\def\xxx{\kern1.7pt}\def\yyy{\kern-1.7pt}\luz{#1}}
\def\spa#1{\def\xxx{\kern1pt}\def\yyy{\kern-1pt}\luz{#1}}

\newcount\beg
\newbox\aabox
\newbox\atbox
\newbox\fpbox
\def\abbrevauthors#1{\setbox\aabox=\hbox{\sevenrm\uppercase{#1}}}
\def\abbrevtitle#1{\setbox\atbox=\hbox{\sevenrm\uppercase{#1}}}
\long\def\pag{\beg=\pageno
\def\leftheadline{\noindent\rlap{\nine\folio}\hfil\copy\aabox\hfil}
\def\rightheadline{\noindent\hfill\copy\atbox\hfill\llap{\nine\folio}}
\def\phead{\setbox\fpbox=\hbox{\sevenrm 
************************************************}%
\noindent\vbox{\sevenrm\baselineskip9pt\hsize\wd\fpbox%
\centerline{***********************************************}

\centerline{BANACH CENTER PUBLICATIONS, VOLUME **}

\centerline{INSTITUTE OF MATHEMATICS}

\centerline{POLISH ACADEMY OF SCIENCES}

\centerline{WARSZAWA 19**}}\hfill}
\footline{\ifnum\beg=\pageno \hfill\nine[\folio]\hfill\fi}
\headline{\ifnum\beg=\pageno\phead
\else
\ifodd\pageno\rightheadline \else \leftheadline \fi 
\fi}}

\newbox\tbox
\newbox\aubox
\newbox\adbox
\newbox\mathbox

\def\title#1{\setbox\tbox=\hbox{\let\\=\cr 
\baselineskip14pt\vbox{\twbf\tabskip 0pt plus15cc
\halign to\hsize{\hfil\ignorespaces \uppercase{##}\hfil\cr#1\cr}}}}

\newbox\abbox
\setbox\abbox=\vbox{\vglue18pt}

\def\author#1{\setbox\aubox=\hbox{\let\\=\cr 
\nine\baselineskip12pt\vbox{\tabskip 0pt plus15cc
\halign to\hsize{\hfil\ignorespaces \uppercase{\spa{##}}\hfil\cr#1\cr}}}%
\global\setbox\abbox=\vbox{\unvbox\abbox\box\aubox\vskip8pt}}

\def\address#1{\setbox\adbox=\hbox{\let\\=\cr 
\nine\baselineskip12pt\vbox{\it\tabskip 0pt plus15cc
\halign to\hsize{\hfil\ignorespaces {##}\hfil\cr#1\cr}}}%
\global\setbox\abbox=\vbox{\unvbox\abbox\box\adbox\vskip16pt}}

\def\mathclass#1{\setbox\mathbox=\hbox{\footnote{}{1991 
{\it Mathematics Subject Classification}\/: #1}}}

\long\def\maketitlebcp{\pag\unhbox\mathbox
\footnote{}{The paper is in final form and no version 
of it will be published elsewhere.} 
\vglue7cc
\box\tbox
\box\abbox
\vskip8pt}

\long\def\abstract#1{{\nine{\bf Abstract.} 
#1

}}

\def\section#1{\vskip-\lastskip\vskip12pt plus2pt minus2pt
{\bf #1}}

\long\def\th#1#2#3{\vskip-\lastskip\vskip4pt plus2pt
{\sc #1} #2\hskip-\lastskip\ {\it #3}\vskip-\lastskip\vskip4pt plus2pt}

\long\def\defin#1#2{\vskip-\lastskip\vskip4pt plus2pt
{\sc #1} #2 \vskip-\lastskip\vskip4pt plus2pt}

\def\Proof{\vskip-\lastskip\vskip4pt plus2pt 
\sp{Proo{f.}\ }\ignorespaces}

\def\endproof{\nobreak\kern5pt\nobreak\vrule height4pt width4pt depth0pt
\vskip4pt plus2pt}

\newbox\refbox
\newdimen\refwidth
\long\def\references#1#2{{\nine
\setbox\refbox=\hbox{\nine[#1]}\refwidth\wd\refbox\advance\refwidth by 12pt%
\def\textindent##1{\indent\llap{##1\hskip12pt}\ignorespaces}
\vskip24pt plus4pt minus4pt
\centerline{\bf References}
\vskip12pt plus2pt minus2pt
\parindent=\refwidth
#2

}}

\def\footnoterule{\kern -3pt \hrule width 4cc \kern 2.6pt}

\catcode`@=11
\def\vfootnote#1%
{\insert\footins\bgroup\nine\interlinepenalty\interfootnotelinepenalty%
\splittopskip\ht\strutbox\splitmaxdepth\dp\strutbox\floatingpenalty\@MM%
\leftskip\z@skip\rightskip\z@skip\spaceskip\z@skip\xspaceskip\z@skip%
\textindent{#1}\footstrut\futurelet\next\fo@t}
\catcode`@=12

\hoffset 1cm  
\def\mbox{\hbox}
\def\ttimes{\mbox{$\hskip.5mm\bigcirc\hskip-4.05mm
\perp\hskip1mm$}}
\def\Ttimes{\mbox{$\hskip.5mm\bigcirc\hskip-3mm
\raise -.5mm \hbox{$\top$}\hskip1mm$}}

\def\id{\mbox{\rm id}\, }
\def\End{\mbox{\bf End}\, }


\font \msa=msam10 scaled \magstep1
\font \msb=msbm10 scaled \magstep0

\def\bC{\mbox{\msb C} }

\def\isp{\mbox{\msa t} }

\font \eul=eufm10 scaled \magstep2

\def\gotG{\mbox{\eul g}}

\def\sfR{\mbox{\bf R} }

\def\Dr{\Delta }
\def\Drr{{\widetilde{\Delta}} }
\def\Drop{\Delta ^{\rm op}}
\def\tR{{\widetilde{R}} }
\def\tE{{\widetilde{E}} }
\def\er{\varepsilon }

\mathclass{Primary 17B37, 16W30; Secondary 81S05.}

\abbrevauthors{S. Zakrzewski}
\abbrevtitle{Coboundary Hopf algebras}

\title{A characterization of coboundary\\
 Poisson Lie groups and Hopf algebras}

\author{Stanis{\L}aw\ Zakrzewski}
\address{Department of Mathematical Methods in Physics,
University of Warsaw\\
Ho\.{z}a 74, 00-682 Warsaw, Poland\\
E-mail: szakrz@fuw.edu.pl}

\maketitlebcp

\footnote{}{Research supported by KBN grant 2 P301 020 07.}

\abstract{We show that a Poisson Lie group $(G,\pi )$ is
coboundary if and only if the natural action of $G\times G$ on $M=G$
is a Poisson action for an appropriate Poisson structure on $M$ (the
structure turns out to be the well known $\pi _+$). We analyze
the same condition in the context of Hopf algebras.  Quantum
analogue of the $\pi _+$ structure on $SU(N)$ is described in
terms of generators and relations as an example.}

\section{1. Preliminaries.}
For the theory of Poisson Lie groups we refer to [1, 2, 3, 4, 5]. 
We follow the notation used in our previous papers [6, 7].

A {\it Poisson Lie group} is a Lie group $G$ equipped with a
Poisson structure $\pi $ such that the multiplication map is
Poisson. The latter property is equivalent to the following
property (called {\it multiplicativity} of $\pi$)
$$
\pi (gh)=\pi (g)h+g\pi (h) \qquad \mbox{for}\;\; g,h\in G.\leqno
(1)
$$
Here $\pi (g)h$ denotes the right translation of $\pi (g)$ by
$h$ etc. This notation will be used throughout the paper.

A Poisson Lie group is said to be {\it coboundary} if 
$$
\pi (g) = rg -gr\leqno
(2)
$$
for a certain element $r\in  \gotG\bigwedge \gotG $. Here $\gotG
$ denotes the Lie algebra of $G$. Any bivector field of the form (2) 
is multiplicative. It is Poisson if and only if
$$
[r,r]\in (\gotG\bigwedge\gotG\bigwedge\gotG )_{\rm inv}
$$
(the Schouten bracket $[r,r]$ is $\gotG $-invariant). 
In this case the element $r$ is said to be a {\it classical
$r$-matrix} ({\it on} $\gotG $).

For any Poisson Lie group $(G,\pi )$, the antipode map $g\mapsto
Sg:=g^{-1}$ is anti-Poisson: 
$$
S_* \pi = -\pi .\leqno
(3)
$$

\section{2. Gauge transformations of a lattice connection on one
link.} 
Consider the following action
$$
(G\times G)\times G \ni ((g_0,g_1),g)\mapsto g_1gg_0^{-1} \in G\leqno
(4)
$$
of $G\times G$ on $G$. This type of action is familiar in gauge
field theory on the lattice. We think here about an `elementary'
lattice composed of only one link with two ends: $0$ and $1$.
Elements $g_0$ and $g_1$ are the values of the gauge
transformation at the lattice sites $0$ and $1$, respectively. The
connection on the link is represented by the element $g$.

One can ask if it is possible to consider the gauge group to be 
a Poisson Lie group (or, a quantum group). In this case it is
natural to require  the action (4) to be a Poisson action
(i.e. the map (4) to be a Poisson map). 

\defin{Definition}{1. 
A Poisson Lie group $(G,\pi )$ is said to be {\it gauge-admissible}\/
 if there exists a Poisson structure $\rho $ on $G$
such that the map (4) is a Poisson map as a map from
$(G,\pi )\times (G,\pi )\times (G,\rho )$ to $(G,\rho)$.}

Note that we treat the gauge group differently than the space of
connections (even if the latter is parameterized by the group
manifold).

\th{Proposition}{1.}{A Poisson Lie group is  gauge admissible if
and only if it is coboundary.}

\Proof
Let $(G,\pi )$ be a Poisson Lie group. It is gauge admissible if
and only if the map
$$
 G\times G\times G \ni (x,y,z)\mapsto xyz^{-1}\in G\leqno
(5)
$$
is Poisson as a map from $(G,\pi)\times (G,\rho )\times (G,\pi
)$ to $(G,\rho )$ or, equivalently (using (3)), if the map
$\Psi \colon G\times G\times G\to G$ defined by
$$
\Psi (x,y,z) = xyz
$$
is Poisson as a map from $(G,\pi)\times (G,\rho )\times (G,-\pi
)$ to $(G,\rho )$. By a similar reasoning which leads to
(1), this is equivalent to
$$
\rho (xyz)= \pi (x)yz + x\rho (y)z -xy\pi
(z)\qquad\qquad\mbox{for}\;\;\; x,y,z\in G.\leqno
(6)
$$
We have two following particular cases of this equality. If we
set $z=e$ (the group unit), we get
$$
\rho (xy)= \pi (x) y + x\rho (y), \leqno
(7)
$$
and if we set $x=e$, we get
$$
\rho (yz) = \rho (y) z-y\pi (z).\leqno
(8)
$$
It is easy to see that (7) and (8) together are equivalent
to (6).  Since $\rho =\pi $ is a particular solution of
(7), the general solution of (7) is given by
$$
\rho (g) = \pi (g) + gA\, ,\leqno
(9)
$$
where $A\in  \gotG\bigwedge\gotG $. Since $\rho =-\pi $ is a
particular solution of 
(8), the general solution of (8) is given by
$$
\rho (g) = -\pi (g) + Bg\, ,\leqno
(10)
$$
where $B\in  \gotG\bigwedge\gotG $. For the compatibility of (9)
and (10) we must have 
$$
\pi (g) ={Bg - gA \over 2}.
$$
Since $\pi (e)=0$, we have $B=A$, and finally
$$
\pi (g) ={Ag - gA\over 2}, \qquad  \qquad \rho (g) ={Ag+gA\over 2}.
$$    
This shows that $(G,\pi )$ is gauge-admissible if and only if
it is coboundary (with $r=A/2$; note that if $r$ is the
classical $r$-matrix then $\pi _+(g):=rg+gr=\rho (g)$ is
automatically a Poisson bivector field).
\endproof

It is clear that for a given coboundary Poisson structure $\pi$, 
all possible $\rho $ are obtained from one by adding an
invariant element of $\gotG\bigwedge\gotG $. In particular, if
$\gotG $ is
semisimple, then $\rho$ is unique.

\section{3. Hopf algebra case.}
 
 Let $(H,m,\Dr )$ be a Hopf algebra. Here $m\colon H\otimes H\to
H$ and $\Dr\colon H\to H\otimes H$ denote the multiplication and the
comultiplication in $H$. Let $I$ and $c$ denote the unit and
counit of the Hopf algebra.

We set
$$
\Psi : = m(m\otimes \id)=m(\id \otimes m)
$$
and ask when there exists a (new) coalgebra structure $\Drr$
(with the same counit $c$) on
$H$ such that $\Psi $ is a morphism from $(H,\Dr )\otimes (H,\Drr
)\otimes (H,\Drop )$ to $(H,\Drr )$. Here $\Drop$ is the
comultiplication opposite to $\Dr $: $\Drop = P\circ\Dr$, where
$P$ is the permutation in the tensor product.

The condition for $\Psi$ to be such a morphism reads:
$$
\Drr \Psi = (\Psi\otimes \Psi) (\id\otimes\id\otimes P\otimes
\id\otimes\id) (\id\otimes P\otimes P\otimes \id )
(\Dr\otimes\Drr\otimes \Drop ),\leqno
(11)
$$
and is equivalent to two following conditions
$$
\Drr m = (m\otimes m)(\id \otimes P\otimes \id )
(\Dr\otimes \Drr )\leqno
(12)
$$
$$
\Drr m = (m\otimes m)(\id \otimes P\otimes \id )
(\Drr\otimes \Drop )\leqno
(13)
$$
(they follow from (11) by applying it to
$\id\otimes\id\otimes I$ and $I\otimes \id\otimes\id $,
respectively).
It is easy to solve these conditions for $\Drr$. 
Applying (12) to $\id\otimes I$, we get 
$$
\Drr (a) = \Dr (a)\, \sfR \qquad\qquad a\in H\, ,\leqno
(14)
$$
where the multiplication is that of $H\otimes H$ and
$$
 \sfR := \Drr (I).
$$
It is easy to see that  (14) solves (12) for any $\sfR $.

Similarly, applying (13) to $I\otimes \id $, we get 
$$
\Drr (a) = \sfR \,\Drop (a) \qquad\qquad a\in H.
$$
This is a solution of (13) for any $\sfR $. It follows that
 the general solution of (11) is
(14), where the $R $-{\it matrix} $\sfR $ satisfies the
compatibility condition
$$
\Dr (a)\, \sfR  = \sfR \, \Drop (a) \qquad\qquad a\in H.\leqno
(15)
$$
It is easy to see that $\Drr $ is coassociative if and only if 
$$
[(\Dr \otimes \id )\,\sfR \,](\,\sfR \otimes I)=[(\id\otimes \Dr )\,\sfR
\, ](I\otimes \sfR ). \leqno
(16)
$$
Indeed,
$$ (\Drr\otimes \id )\Drr (a)= [(\Dr \otimes \id )(\Dr (a)
\sfR )](\sfR \otimes \id )= [(\Dr \otimes \id )\Dr (a)][(\Dr \otimes \id
)\sfR ](\sfR \otimes \id ),$$
$$ (\id \otimes\Drr )\Drr (a)=
[(\id\otimes \Dr )(\Dr (a)\sfR )](\id\otimes \sfR )=
[(\id\otimes \Dr )\Dr (a)][(\id\otimes \Dr )\sfR ](\id\otimes \sfR ).$$
Concluding: the question at the beginning of this section has an
affirmative answer if and only if there exists an element $\sfR \in
H\otimes H$ such that (15), (16) hold and
$$
(c\otimes\id )\, \sfR  = I = (\id\otimes c)\, \sfR .
$$
A Hopf algebra satisfying those conditions might be called 
{\it gauge-admissible}\/ or {\it coboundary}\/. (I do not know whether
one can always choose $\sfR $ to be `unitary', like in [2]:
$\sfR _{12}\sfR _{21}= I\otimes I$).

The Hopf algebra considered in this section should be
interpreted as a dual of the Hopf algebra of functions on a
quantum group (quantized universal enveloping algebra). In the
next section we give an example of a `gauge-admissible' matrix
quantum group.  

\section{4. Example in terms of generators and relations.}

Let 
$$
R (u\Ttimes u) = (u\Ttimes u)R\leqno
(17)
$$
be a part of relations defining a matrix quantum group $(A,u)$. Here
$u= (u_{ij})_{i,j=1,\ldots ,n}$ is the defining representation of
the quantum group, $R$ is the fundamental intertwiner
($R$-matrix of FRT-type) and we use the Woronowicz's notation for the  
`matrix' tensor product. Let us note that we have
$$
\tR (u ^{-1}\Ttimes u^{-1}) = (u ^{-1}\Ttimes u^{-1}) \tR\, ,\leqno
(18)
$$
where $\tR := PRP$. Let us denote by $B$ the algebra generated
by the entries of the $n\times n$ matrix $w$ and relations
$$
R (w\Ttimes w) = (w\Ttimes w)\tR.\leqno
(19)
$$
It is easy to see that there exists exactly one homomorphism $\isp$
(quantum gauge transformation -- the analogue of (5))
 from $B$ to $A\otimes B\otimes A$ such that
$$
\isp (w^i{_j})= \sum_{kl}u^i{_k}\otimes w^k{_l}\otimes
(u^{-1})^l{_j}\, , 
$$
or, using the Woronowicz's notation,
$$
(\isp\otimes \id) (w) = u\ttimes w\ttimes u^{-1}
$$
(here $w$ is understood as an element of $\End (\bC ^n)\otimes B$).
In order to see that $u\ttimes w\ttimes u^{-1}$ satisfies the
same relations as $w$, we notice that
$$ (u\ttimes w\ttimes u^{-1})\Ttimes (u\ttimes w\ttimes u^{-1})=
(u\Ttimes u)\ttimes (w\Ttimes w)\ttimes (u^{-1}\Ttimes u^{-1})$$
and use subsequently (17), (18) and (19).

In order to be more precise, we consider now a specific matrix
quantum group, namely $SU_q(n)$, as given in [8].
The *-algebra $A$ of `regular functions' on $SU_q(n)$ is the one generated
by the entries of an $n\times n$ matrix $u$ and the following relations:
$$
u^{(n)}E=E,\qquad E'u^{(n)}=E', \qquad uu^* = I_n\otimes I_A = u^*u .
\leqno
(20)
$$
Here $u^{(n)}$ is the $n$-th tensor power of $u$,
$E$ is the `$q$-deformed'
volume element
$$
      E^{i_1i_2\ldots i_n}= 
(-q)^{ {\rm number\;\, of\;\, inversions\;\, in }\;\; (i_1,\ldots,i_n)}
,\qquad\qquad E'_{i_1\ldots i_n}= E^{i_1\ldots i_n}
$$
(for $(i_1\ldots i_n)$ not being a permutation we set
$E^{i_1\ldots i_n}=0$) and
$I_n$ is the unit $n\times n$ matrix. Note that in this case
$$
 (u^{-1})^{(n)}\tE = tE, \qquad \tE ' (u^{-1})^{(n)}= \tE ',
$$
where
$$
\tE = P_{\rm total}\, E,\qquad \tE ' = E' P_{\rm total},
$$
 $P_{\rm total}$ being the total permutation $(1,2,\ldots ,n)\mapsto
(n,\ldots ,2,1)$.

Let $B$ be the *-algebra generated by the entries of an $n\times
n$ matrix $w$ and relations
$$
w^{(n)}\tE=(-1)^{{n(n-1)\over 2}}E,\qquad
E'w^{(n)}=(-1)^{{n(n-1)\over 2}}\tE ', \qquad ww^* = I_n\otimes
I_B = w^*w . \leqno
(21)
$$
It is easy to check that $u\ttimes w\ttimes u^{-1}$ satisfies
the same relations, hence we have the `gauge transformations' on
the quantum level.

It is essential to know if algebra $B$ has a correct size
(Poincar\'{e} series), i.e. if the deformation is {\it flat}. We
shall show that $B$ is actually isomorphic to $A$. To this end,
consider the change of variables
$$ 
   u =\er w P_{\rm total} 
$$
in (20), where $\er $ is such a complex number that $\er
^n = (-1)^{{n(n-1)\over 2}}=\det P_{\rm total}$. It is easy to
see that relations (20) are now transformed to relations
(21). 

\section{5. Remarks.}
\vskip4pt plus2pt

{\bf 5.1.} \  
The algebra $B$ defined in (21) is the quantum
counterpart of the Poisson structure $\pi _+ (g)=rg+gr$ on $SU(n)$.
The case of a general group is sketched in (19). Note that
if we substitute $u=wg_0$ in (17) where $g_0$ is an element
of the classical group such that 
$$
(g_0\otimes g_0)R(g_0^{-1}\otimes g_0^{-1})= PRP\, ,\leqno
(22)
$$
then we obtain relations (19). One can check that the
well known $R$-matrix for the $A_n$ series satisfies (22)
if we choose $g_0=\er P_{\rm total}$. The corresponding fact for
Poisson groups means that we find an element $g_0\in G$ such that
$$
\pi = \pi _+ g_0\, \leqno
(23)
$$
 i.e. $\pi (gg_0) = \pi _+ (g)g_0$, that is to say
$$
rgg_0-gg_0r = rgg_0+grg_0,
$$
or,
$$
       g_0rg_0^{-1} = -r,
$$
 or,
$$
    \pi _+ (g_0)=0.\leqno
(24)
$$
For instance in the case of the standard $r$-matrix of the $A_n$-series,
$$ r = \sum_{j<k} e_j{^k}\wedge e_k{^j},$$
$g_0:=\er P_{\rm total}$ will do the job, because $Pe_j=e_{j'}$,
$j':=n+1-j$.
\vskip4pt plus2pt

{\bf 5.2.} \  Formula (14) was used in [9] to discuss
twisting Hopf algebras by 2-cocycles. The Poisson structure $\pi
_+$ is isomorphic to $\pi$ by a translation (23) if
and only if it vanishes at some point (namely $g_0$, see (24)).
This situation (and previously discussed isomorphism of $B$ with
$A$) corresponds to twisting by a coboundary.
\vskip4pt plus2pt

{\bf 5.3.} \ The author would like to thank to Jiang-Hua Lu, Shahn
Majid and Marco Tarlini for enlightening discussions.

\references{9}{
\item{[1]} V. G. Drinfeld, {\it Hamiltonian structures on 
Lie groups, Lie bialgebras and the meaning of the classical
Yang-Baxter equations}\/,  Soviet~Math.~Dokl.~27 (1983),
68--71.

\item{[2]} V. G. Drinfeld, {\it Quantum groups}\/, Proc. ICM, 
Berkeley, 1986, vol.1, 789--820.

\item{[3]} M. A. Semenov-Tian-Shansky, {\it Dressing
transformations and Poisson Lie group actions}\/,
Publ. Res. Inst. Math. Sci.,
Kyoto University 21 (1985), 1237--1260.

\item{[4]} J.-H. Lu and A. Weinstein, {\it Poisson Lie Groups,
Dressing Transformations and Bruhat Decompositions}\/,
J.~Diff.~Geom.~31 (1990), 501--526.

\item{[5]} J.-H. Lu, {\it Multiplicative and affine Poisson
structures on Lie groups}\/, Ph.D. Thesis, University of
California, Berkeley (1990).

\item{[6]} S. Zakrzewski,
{\it Poisson structures on the Lorentz group}\/, 
 Lett.~Math.~Phys.~32 (1994), 11--23.

\item{[7]} S.~Zakrzewski, {\it Poisson homogeneous spaces}\/,
in: ``Quantum Groups, Formalism and Applications'', 
Proceedings of the XXX Winter School on Theoretical Physics
14--26 February 1994, Karpacz, J.~Lukierski, Z.~Popowicz,
J.~Sobczyk (eds.), Polish Scientific Publishers PWN, Warsaw
1995, pp.~629--639.

\item{[8]} S.L. Woronowicz, {\it Tannaka-Krein duality for
compact matrix pseudogroups. Twisted SU(N) groups}\/,
Invent.~Math.~93 (1988), 35.

\item{[9]} J.-H. Lu, {\it On the Drinfeld double and the
Heisenberg double of a Hopf algebra}\/, 
Duke Math. Journ. 74, No.3 (1994), 763--776.}

\bye